\documentclass{svmult}

\usepackage{multicol}        
\usepackage[bottom]{footmisc}

\def\pt{\forall\;}
\def\RR{{\cal R}}
\def \th{\mbox{th}\,}
\def \ch{\mbox{ch}\,}

\def \ddt#1{\frac{\partial #1}{\partial t}}

\def \dd#1#2{\frac{\partial #1}{\partial #2}}
\def \ddx#1{\frac{\partial #1}{\partial x}}
\def\d2x#1{\frac{\partial ^2 #1}{\partial x^2}}      
\def\div{\mbox{{\rm div}}\;}
\def\rot{\mbox{{\rm curl}}\;}
\def\ds{\longrightarrow}
\def\E{{\mathcal E}}

\def\L{{\mathcal L}}

\def\R{{\mathord{I\!\! R}}}

\def\d{\delta}

\def\M{{\mathcal M}}
\def\K{{\mathcal K}}
\def\v{\wedge}

\def\e{\varepsilon} 
\def\dsp{\displaystyle}

\begin{document}

\title*{Stability for Walls in Ferromagnetic Nanowire}

\author{G. Carbou \inst{1}\and
S. Labb\'e \inst{2}}

\institute{MAB, Universit\'e Bordeaux 1, 351 cours de la Lib\'eration,
33405 Talence cedex
\texttt{carbou@math.u-bordeaux1.fr}
\and Universit\'e Paris Sud, Orsay \texttt{labbe@math.u-psud.fr}}

\maketitle

\begin{abstract}
We study the stability of travelling wall profiles for a one
dimensional model of ferromagnetic nanowire submitted to an exterior
magnetic field. We prove that
these profiles are asymptotically stable modulo a translation-rotation
for small applied magnetic fields.
\end{abstract}

\section{Model for ferromagnetic nanowires}
\label{sec:1}

Ferromagnetic materials are characterized by a spontaneous
magnetization described by the magnetic moment $u$ which is a unitary
vector field linking the magnetic induction $B$ with the magnetic
field $H$ by the relation $B=H+u$. The variations of $u$ are described by the Landau-Lifschitz Equation
\begin{equation}
\label{LL3D}
\frac{\partial u}{\partial t}=-u\wedge H_e-u\wedge (u\wedge H_e)
\end{equation}
where the effective field is given by $H_e=\Delta u + h_d(u)+ H_a$, 
and  the demagnetizing field $h_d(u)$  is deduced from $u$ solving the
magnetostatic equations:
\[\div B=\div(H+u)=0 \mbox{ and } \rot H=0\]
where $H_a$ is an appplied magnetic field. 

For more details on the
 ferromagnetism model, see \cite{brown}, \cite{Halpern.Labbe:Modelisation}, \cite{LL} and
\cite{W}. For existence results about Landau Lifschitz equations see
\cite{carbou.fabrie:time}, \cite{carbou.fabrie:regular}, \cite{JMR},
\cite{V}. For numerical studies see \cite{Had}, \cite{L} and
\cite{Labbe.Bertin:Microwave}. For asymptotic studies see
\cite{Alouges:Neel}, \cite{CFG}, \cite{desimone.kohn.muller.otto:magnetic},
\cite{riviere.serfaty:compactness} and  \cite{sanchez}.

\vspace{4mm}

In this paper we consider an asymptotic one dimensional model of ferromagnetic
nanowire submited to an applied field along the axis of the wire.  We
denote by  $(e_1,e_2,e_3)$  the canonical basis of $\R^3$.  The
ferromagnetic nanowire is assimilated to the axis $\R e_1$. The demagnetizing energy is
approximated by the formula $h_d(u)=-u_2 e_2 - u_3 e_3$ where $u=(u_1,u_2,u_3)$ (this
approximation of the demagnetizing energy for a ferromagnetic wire is
obtained using a BKW method by D. Sanchez, taking the limit when the
diameter of the wire tends to zero in \cite{sanchez}). We assume in
addition that an exterior magnetic field $\delta e_1$ is applied
along the wire axis. 

To sum up  we study the following system 
\begin{equation}
\label{LL1D}
\left\{
\begin{array}{l}
\dsp \frac{\partial u}{\partial t}=-u\wedge h_\d (u) -u\wedge(u\wedge h_\d(u))\\
\mbox{with } h_\d(u)=\dsp \dd{^2u}{x^2}-u_2 e_2 - u_3 e_3 + \delta e_1
\end{array}\right.
\end{equation}

For $\delta=0$, that is without applied field, we observe in physical
experiments the formation of a wall breaking down the domain in two parts:
one in which the magnetization is almost equal to $e_1$ and another in which
the magnetization is almost equal to $-e_1$. Such a distribution is
described in our one dimensional model by the following profile $M_0$:

\begin{equation}
M_0=
\left(
\begin{array}{c}
\th{x}\\ 0 \\ \frac{1}{\ch x}\\
\end{array}
\right).
\end{equation}

This profile is a steady state solution of Equation (\ref{LL1D}) with
$\delta=0$. We prove in \cite{CL} the stability of the profile $M_0$
for Equation (\ref{LL1D}) without applied field (when $\delta=0$). 

\vspace{2mm}
When we apply a magnetic field in the direction $+ e_1$ (that is with
$\delta >0$) since Landau-Lifschitz Equation tends to align the
magnetic moment with the effective field,   we
observe a translation of the wall in the direction $- e_1$. Furthermore, we 
observe a rotation of the magnetic moment around the wire axis. This phenomenon is described by the
solution of (\ref{LL1D})

\begin{equation}
\label{Udelta}
U_\delta(t,x)=R_{\delta t}(M_0(x+\delta t))
\end{equation}

where $R_\theta$ is the rotation of angle $\theta$ around the axis $\R
e_1$:
\[R_\theta=\left(
\begin{array}{ccc}
1&0&0\\
\\
0&\cos \theta &-\sin \theta\\
\\
0&\sin\theta & \cos \theta
\end{array}
\right)
\]

We study in this paper the stability of $U_\delta$,
we prove that for a small $\delta$, $U_\delta$ is stable for
the $H^2$ norm and asympotically stable for the $H^1$ norm, modulo a
translation in the variable $x$ and a rotation around $\R e_1$. This
result is claimed in the following theorem:

\begin{theorem}
\label{thm}
There exists $\delta_0>0$ such that for all $\delta$ with
$|\delta|<\delta_0$ then for $\e>0$ there exists $\eta>0$ such that if
 $\|u(t=0,x)-U_\delta(t=0,x)\|_{H^2}<\eta$ then the solution $u$ of
Equation (\ref{LL1D}) with initial data $u(t=0,x)$ satisfies:
\[\pt t>0, \|u(t,x)-U_\delta(t,x)\|_{H^2}<\e.\]
In addition there exists
$\sigma_\infty$ and $\theta_\infty$ such that
\[\|u(t,x)-R_{\theta_\infty}(U_\delta(t,x+\sigma_\infty))\|_{H^1}\ds
0\mbox{ when } t \ds +\infty.\]
\end{theorem}

This result is a generalization of the stability result concerning the
static walls when $\delta=0$ in \cite{CL}. It looks like the theorems of stability concerning the
travelling waves solutions for semilinear equations like Ginzburg
Landau Equation (see Kapitula \cite{kapitula:multidimensionnal}). Here
we have three new difficulties. The first one is that the magnetic
moment takes its values in the sphere and not in a linear space. In
order to work with maps with values in a linear space we will use a
mobile frame adapted to Landau-Lifschitz equation and we will describe
in Section 2 the magnetic moment in this mobile frame. The second difficulty is
that we have here a two dimensional invariance family for Equation
(\ref{LL1D}) whereas Ginzburg Landau-Equation is only invariant by
translation. This is the reason why we must use in the perturbations
description the translations and the rotations (see Section 3). The last difficulty is
that Landau-Lifschitz Equation is quasilinear, and then we have to
couple variational estimates and semi-group estimates to control the
perturbations of our profiles. Section 4 is devoted to these estimates.

\section{Landau-Lifschitz Equation in the mobile frame}

\subsection{First reduction of the problem}

For $u$ a solution of Landau-Lifschitz Equation (\ref{LL1D}) we define
$v$ by $v(t,x)=R_{-\delta t}(u(t,x-\delta t))$ (that is $u(t,x)=R_{\delta t}(v(t,x+\delta t))$). A straightforward
calculation gives that $u$ satisfies (\ref{LL1D}) if and only if $v$ satisfies 
\begin{equation}
\label{LL2}
\left\{
\begin{array}{l}
\dsp 
\ddt{v}=-v\v h(v)-v\v(v\v h(v))-\delta(\dd{v}{x}+ v_1 v - e_1)
\\ 
\dsp h(v)=\dd{^2v}{x^2}- v_2 e_2 - v_3 e_3
\end{array}\right.
\end{equation}

In addition  $U_\delta$ is stable for (\ref{LL1D}) if and only if $M_0$ is stable for
(\ref{LL2}), that is we are led to study the stability of a static
profile, which is more convenient.

\subsection{Mobile frame}

Let us introduce the   mobile frame $(M_0(x), M_1(x), M_2)$, where 
\[M_1(x)=\left(\begin{array}{c}
\dsp \frac{1}{\ch x}\\0\\ -\th
x\end{array}\right)\mbox{ and
}M_2=\left(\begin{array}{c}0\\1\\0\end{array}\right)\]
Let $v:\R^+_t\times \R_x\ds S^2\subset \R^3$ be  a little perturbation of $M_0$. We can decompose $v$ in the
mobile frame writting 
\[v(t,x)=r_1(t,x) M_1(x) + r_2(t,x) M_2 + \sqrt{1-r_1^2-r_2^2}M_0(x).\]

Now we can obtain a new version of Landau-Lifschitz Equation: $v$ satisfies (\ref{LL2}) if and only if 
 $r=(r_1,r_2)$ satisfies
\begin{equation}
\label{LL3}
\ddt{r}=(\L+ \delta l) r+G(r)(\dd{^2r}{x^2})+ H(x,r,\dd{r}{x})
\end{equation}

where
\begin{itemize}
\item the linear operator $\L$ is given by  $\L=J L$ with  $\dsp J=\left(\begin{array}{cc}
-1&-1\\1&-1\end{array}\right)$
and  
\[\dsp L=-\dd{^2}{x^2} +2\th ^2 x-1\],

\item the linear perturbation due to the presence of the applied magnetic
 field $\delta e_1$ is given by 
$\delta l$ with  $\dsp l=\dd{}{x}+\th x$,

\item the higher degree non linear part is $G(r)(\dd{^2r}{x^2})$,
where $G(r)$ is a matrix depending on $r$ with $G(0)=0$,

\item the last non linear term $H(x,r,\dd{r}{x})$ is at least
quadratic in the variable $\dsp (r,\dd{r}{x})$.

\end{itemize}

In addition the stability of the profile $M_0$ for Equation
(\ref{LL2}) is equivalent to the stability of the zero solution for
Equation (\ref{LL3}).

\section{A new system of coordinates}

We remark that $L$ is a self adjoint operator on $L^2(\R)$, with
domain $H^2(\R)$. Furthermore, $L$ is positive since we can write
$L=l^*\circ l$ with $l=\dsp \ddx{}+\th x$, and Ker $L$ is the one
dimensional space generated by $\dsp \frac{1}{\ch x}$.

The matrix $J$ being invertible, Ker $\L$ is the two dimensional
space generated by $v_1$ and $v_2$ with 
\[
v_1(x)=\left(\begin{array}{c}
0\\\dsp\frac{1}{\ch x}\end{array}\right)
,\;\;\;
v_2(x)=\left(\begin{array}{c}
\dsp\frac{1}{\ch x}\\ 0\end{array}\right)
\]

We introduce $\E=($Ker $\L )^\perp$. We denote by $Q$ the orthogonal
projection onto $\E$ for the $L^2(\R)$ scalar product.

\vspace{3mm}

Landau-Lifschitz equation (\ref{LL2}) is invariant by  translation in the
variable $x$ and by rotation around the axis $e_1$. Therefore for
$\Lambda=(\theta,\sigma)$ fixed in $\R^2$, $M_\Lambda$ defined by
$M_\Lambda(x)=R_\theta(M_0(x-\sigma))$ is a solution of Equation
(\ref{LL2}). We introduce $R_\Lambda(x)$ the coordinates of
$M_\Lambda(x)$ in the mobile frame
$(M_1(x),M_2(x))$:
\[R_\Lambda(x)=\left( 
\begin{array}{l}
 M_\Lambda(x)\cdot M_1(x)
\\
M_\Lambda(x)\cdot M_2
\end{array}
\right)\]

\vspace{2mm}

The map $\Psi$ given by
\[\begin{array}{rcl}
\Psi: \R^2\times \E& \ds & H^2(\R)\\
(\Lambda,W) & \longmapsto & r(x)=R_\Lambda (x)+W(x)
\end{array}
\]
is a diffeomorphism in a neighborhood of zero. Thus we can write the
solution $r$ of Equation (\ref{LL3}) on the form :
\[r(t,x)=R_{\Lambda(t)}(x) + W(t,x)\]
where for all $t$, $W(t)\in \E$ and where $\Lambda:\R^+_t\mapsto \R^2$.

\vspace{2mm}
We will re-write Equation (\ref{LL2}) in the coordinates
$(\Lambda,W)$. Taking the scalar product of (\ref{LL2}) with $v_1$ and
$v_2$ we obtain the equation satisfied by $\Lambda$, and using  $Q$
the orthogonal projection onto $\E$, we deduce the equation satisfied
by $W$. After this calculation we obtain that $r$ is solution of
Equation (\ref{LL2}) if and only if $(\Lambda,W)$ satisfies the
following system
\begin{equation}
\label{merde}
\left\{
\begin{array}{l}
\dsp \ddt{W}=(\L + \delta l + \K_{\Lambda}) W + \RR_1(x,\Lambda,W)(\dd{^2W}{x^2})+ 
 \RR_2(x,\Lambda,W,\ddx{W}) 
\\
\\
\dsp
\frac{d\Lambda}{dt}=\M(W,\ddx{W},\Lambda)
\end{array}\right.
\end{equation}
where
\begin{itemize}
\item $\K_\Lambda:H^2(\R)\ds\E$ is a linear map satisfying
\begin{equation}
\label{estiKlambda}
\exists K_1, \;\pt \Lambda\in\R^2, \; \pt W\in\E, \; \|\K_\Lambda
W\|_{L^2(\R)}\leq K_1|\Lambda|\|W\|_{H^2(\R)}
\end{equation}
\item the  non linear terms take their values in $\E$ and satisfy that
there exists a constant $K_2$ such that for $|\Lambda|\leq 1$ and for
all $W\in\E$ 
\begin{equation}
\label{estiR}
\begin{array}{l}
\dsp  \| \RR_1(.,\Lambda,W)(\dd{^2W}{x^2})\|_{L^2(\R)}\leq
K_2\|W\|_{H^1(\R)}\|W\|_{H^2(\R)}\\ \\
\dsp  \| \RR_2(.,\Lambda,W,\ddx{W})\|_{H^1(\R)}\leq
K_2\|W\|_{H^1(\R)}^2
\end{array}
\end{equation}
\item $\M:H^1(\R)\times L^2(\R)\times \R^2\ds \R^2$ satisfies
\begin{equation}
\label{estiM}
\exists K_3, \; \pt \Lambda \mbox{ such that }|\Lambda|\leq 1, \; \pt
W\in \E, \; | \M(W,\ddx{W},\Lambda)|\leq
K_3\|W\|_{H^1(\R)}
\end{equation}

\end{itemize}

Theorem \ref{thm} is equivalent to the following Proposition:

\begin{proposition}
\label{prop}
There exists $\delta_0>0$ such that for $\delta$ with
$|\delta|<\delta_0$, we have the following stability result for
Equation (\ref{merde}): for $\e>0$ there exists $\eta>0$ such that if
$|\Lambda_0|<\eta$ and if $\|W_0\|_{H^2}<\eta$ then the solution
$(\Lambda,W)$ of (\ref{merde}) with initial value $(\Lambda_0,W_0)$ satisfies
\begin{enumerate}
\item 
for all $t>0$, $\|W(t)\|_{H^2}\leq \e$ and $|\Lambda|\leq \e$,

\item $\|W(t)\|_{H^1}$ tends to
zero when $t$ tends to $+\infty$,

\item there exists $\Lambda_\infty\in\R^2$ such that  $\Lambda(t)$
tends to $ \Lambda_\infty$ when $t$ tends to $+\infty$.
\end{enumerate}
\end{proposition}

The last section is devoted to the proof of Proposition \ref{prop}.

\section{Estimates for the perturbations}

\subsection{Linear semi group estimates}

On $\E$ we have
{\it Re} (sp $\L)\subset]-\infty, -1]$. In particular this fact
implies that the $H^2$ norm is equivalent on $\E$ to the norm
$\|\L u\|_{L^2}$. Furthermore it implies good decreasing properties for
the semigroup generated by $\L$. We first prove that this decreasing
property is preserved for the linear part of the Equation on $W$ in
(\ref{merde}) for a little applied field, and if we assume that
$\Lambda$ remains little. 
 
\vspace{2mm}

The operator $l$ is an order one operator dominated on $\E$ by $\L$,
thus there exists $\delta_0>0$ such that if $|\delta|< \delta_0$, {\it
Re} (sp $\L+\delta l)\subset]-\infty, -1/2[$. 

Let us fix $\delta$ such that $|\delta|<\delta_0$. With Estimate
(\ref{estiM}), if $\Lambda$ remains small, $\K_\Lambda$ is a little
perturbation of $\L+\delta l$. This implies that for $\Lambda$ little,
the semigroup generated by $\L+\delta l + Q \K_\Lambda$ has the same
good decreasing properties than $\L$, that is 
there exists $\nu_0>0$ such that if $|\Lambda(t)|$ remains less than $
\nu_0$ for all $t$, then there exists $K_4$ and $\beta>0$ such that
\begin{equation}
\label{esti}
\begin{array}{ll}
\|S_\Lambda (t)W_0\|_{H^1}&\leq K_4\dsp  e^{-\beta t}\|W_0\|_{H^1}
\\
& \leq  K_4\dsp \frac{ e^{-\beta t}}{\sqrt{t}}\|W_0\|_{L^2}.
\end{array}
\end{equation}

We can then use the Duhamel formula to solve the equation on $W$ in
(\ref{merde}):

\[W(t)=S_\Lambda(t) W_0 + \int_0^t S_\Lambda(t-s) \RR_1(s) ds +
\int_0^t S_\Lambda(t-s) \RR_2(s) ds\]
and then using the estimates (\ref{estiR}) and (\ref{esti}) we obtain that if
$|\Lambda(t)|$ remains less than $\nu_0$ then   there exists $K_5$
such that
\begin{equation}
\label{blurp}
\begin{array}{rl}
\dsp \| W(t) \|_{H^1}\leq & \dsp  K_5 e^{-\beta t}\|W_0\|_{H^1} + \int_0^t K_5
\frac{e^{-\beta (t-s)}}{\sqrt{t-s}}\|W(s)\|_{H^1}\|W(s)\|_{H^2}\\ \\
&\dsp +
\int_0^t K_5 e^{-\beta (t-s)}\|W(s)\|_{H^1}^2
\end{array}
\end{equation}

\subsection{Variational estimates}

We see that Estimate (\ref{blurp}) is not sufficient to conclude since
we have the $H^2$ norm of $W$ in the right hand side of this
estimate. In order to dominate this $H^2$ norm, we multiply the
equation on $W$ in (\ref{merde}) by $J
^2 \L^2 W$ and  we obtain that there exists a constant $K_6$:

\[\frac{d}{dt}\|LW\|^2_{L^2} + \|
L^{\frac{3}{2}}W\|^2_{L^2}\left(1-K_6\|LW\|_{L^2}\right)\leq 0\]

>From this estimate we deduce that if $
\|LW\|_{L^2}<\frac{1}{K_6}$, then $1-K_6\|LW\|_{L^2}$ is positive, thus
$ \frac{d}{dt}\|LW\|^2_{L^2}$ is negative and
$\|LW\|_{L^2}$ remains less than $\frac{1}{K_6}$. So if $
\|LW_0\|_{L^2}<\frac{1}{K_6}$, then  for all $t$
$\|LW(t)\|_{L^2}\leq \|LW_0\|_{L^2}$. This property gives a bound for
the $H^2$ norm of $W$ since the $H^2$ norm is equivalent on $\E$ to
$\|LW\|_{L^2}$, and reducing the $H^2$ norm of $W_0$, we obtain the
first part of the conclusion 1 in Proposition \ref{prop}.

\subsection{Conclusion}

Let us assume that $\|L W_0\|_{L^2(\R)}\leq \frac{1}{K_6}$. Then for
all $t$, $\|W(t)\|_{H^2(\R)}\leq C_1 \|LW(t)\|_{L^2}\leq
\|LW_0\|_{L^2}\leq C_2 \|W_0\|_{H^2(\R)}$, where $C_1$ and $C_2$ are constants.

Multiplying  (\ref{blurp}) by $(1+t)^2$, defining $\dsp
G(t)=\max_{[0,T]}(1+s)^2\|W(s)\|_{H^1}$, we obtain that there exists a
constant $K_7$  such that if $|\Lambda(t)|$ remains less than $\nu_0$
we have:
\[G(t)\leq K_7 G(0)+ K_7 G(t) \|W_0\|_{H^2} + K_7 (G(t))^2\]
If we suppose in addition that $\|W_0\|_{H^2}\leq \frac{1}{2K_7}$ we
obtain that
\begin{equation}
\label{pol}
0\leq K_7 G(0)-\frac{1}{2}G(t) + K_7 (G(t))^2\ :=\; P(G(t))
\end{equation}

The polynomial map $P(\xi)=K_7 \xi^2 -\frac{1}{2}\xi + K_7 G(0)$ has
for $G(0)$ small enough two positive roots. We denote by $\xi(G(0))$
the smallest one. For $G(0)$ little enough we have $G(0)\leq
\xi(G(0))\leq 2  K_7 G(0)$ (we can a priori assume that $K_7\geq 1$
for example). Estimate (\ref{pol}) implies that for all
$t$, $G(t)\leq \xi(G(0))$ that is 
\begin{equation}
\label{Wpetit}
\pt t>0, \|W(t)\|_{H^1(\R)}\leq \frac{\xi(G(0))}{1+t^2}\leq \frac{2
K_7 G(0)}{1 + t^2}.
\end{equation}
This implies that $\|W(t)\|_{H^1(\R)}$ tends to zero when $t$ tends to
$+\infty$. It remains to prove that $\Lambda$ remains less that
$\nu_0$ and admits a limit when $t$ tends to $+\infty$.

\vspace{2mm}
Plugging Estimate (\ref{Wpetit}) in 
 the equation on $\Lambda$ in (\ref{merde}) and using (\ref{estiM}),
we obtain that $\frac{d\Lambda}{dt}$ is integrable on $\R^+$, that is
$\Lambda$ admits a limit when $t$ tends to $+\infty$. Furthermore, by
integration we have
\[\pt t, \; |\Lambda(t)|\leq |\Lambda(0)|+\int_0^t K_3
\frac{2K_7 G(0)}{1+s^2}ds\leq |\Lambda(0)|+\pi  K_3 K_7 G(0) \]
Reducing $|\Lambda_0|$ and $G(0)=\|W_0\|_{H^1(\R)}$ we obtain that
for all $t$, $|\Lambda(t)|$ remains less than $\nu_0$, which justifies
all our estimates a posteriori.

\end{document}